\documentclass{article}
\usepackage{amsbsy,amssymb,amscd,amsfonts,latexsym,amstext,delarray,
amsmath}
\newtheorem{thm}{Theorem}[section]

\numberwithin{equation}{section}

\def\ra{\rightarrow}

\def\Q{{\mathbb Q}}

\def\Z{{\mathbb Z}}

\def\d{ \delta }

\def\1ox{{ \Omega^1_{\scriptstyle{X}} }}
\def\2ox{{ \Omega^2_{\scriptstyle{X}} }}

\def\ok1{{ \Omega^1_K }}
\def\ok2{{ \Omega^2_K }}

\def\ra{{ \rightarrow }}
\def\da{{ \downarrow }}

\def\a{{ \alpha }}

\def\e{{ \epsilon }}

\def\hra{{ \hookrightarrow }}
\def\da{{ \downarrow }}

\def\D{{ \Delta }}

\def\8{{ {\infty } }}

\def\^{{ ^{\wedge} }}

\def\Z{{ {\bf Z } }}

\usepackage{amsmath}
\usepackage{amsfonts}
\usepackage{amscd}
\usepackage{amssymb}
\newtheorem{conj}{Conjecture}

\def\bX{{ \bar{X} }}

\def\D{{ \Delta }}

\def\D{{ \Delta }}

\def\bj{ \bar{j} }
\title{The Picard-Lefschetz formula and a conjecture of Kato: The case of
Lefschetz fibrations}
\author{Caterina Consani and Minhyong Kim }

\begin{document}
\maketitle
\section{Introduction}

The rational cohomology groups of  complex algebraic varieties 
 possess
Hodge structures which can be used to refine their
usual topological nature. 
Besides being useful parameters for geometric classification,
these structures reflect in essential ways the
conjectural category of motives, and hence, are of great interest
from the viewpoint of arithmetic. A key link here is provided 
by the Hodge conjecture, which says that for smooth projective
varieties, a class of type $(n,n)$ in $H^{2n}$ should  come from
an algebraic cycle. Stated slightly  more 
abstractly, a copy of the trivial Hodge structure inside
$H^{2n}(X)(n)$
for a smooth projective variety  $X$ should be generated by  an algebraic
cycle class on $X$. Obviously, for the Hodge structures coming from smooth
projective varieties, $H^{2n}(n)$ are the only ones that are of
weight zero, and therefore, that can have trivial sub-structures.

On the other hand, there are many important mixed Hodge structures (MHS)
(\cite{De1} \cite{De2})
that arise from algebraic geometry, including open or singular varieties,
local cohomologies, and the cohomology of Milnor fibers.
The presence of trivial substructures is likely to
reflect subtler phenomena (than the Hodge conjecture!) in the mixed
case, having to do (possibly in a complicated way)
 with trivial substructures of various pure sub-quotients.

The example we have in mind here  is that of limit Hodge structures coming from
degenerations (\cite{St}). That is, let $\D$ be the unit disk and
$$X \ra \D$$
a flat projective family of varieties of fiber dimension $d$
whose only critical value is
the origin and assume that the divisor
$Y:=X_0$ has
 normal crossings. We also assume that the total space is
non-singular of dimension $d+1$
and that the fundamental group of the punctured
base $\pi_1(\D^*) \simeq \Z$ acts unipotently on
the cohomology of a generic fiber.
Let $\D^*=\D-\{0\}$, $X^*=X-X_0$ and
$\bX^*$ be the pull-back of $X^*$ to the universal cover $H$ 
(the upper-half plane)
of
$\D^*$ via the covering map $u\ra \exp (2\pi i u)$. Thus, $\bX^*$ is homotopy
equivalent to any smooth fiber $X_t$ of $X$. Then Steenbrink has put on
$H^n(X_t)\simeq H^n(\bX^*)$ a natural MHS via the theory of nearby cycles.
That is, if we denote by
$\bj: \bX^* \ra X$ the natural map,
$i: Y \hra X$ the inclusion, and the complex of nearby cycles
$$R\Psi:=i^*\bj_* \Q,$$
 one has natural  isomorphisms
$$H^n (X_t, \Q)\simeq H^n(\bX^*, \Q)\simeq H^n(Y,  R\Psi)\simeq
H^n(Y,  (R\Psi)_u)  $$
where $(R\Psi)_u \subset R\Psi$ denotes the maximal subcomplex
on which the monodromy action is unipotent (one can make this
precise
with certain choices of complexes representing the objects in
the derived category).
The complex $(R\Psi)_u $ is shown to underlie a cohomological
mixed Hodge complex (which depends on the choice of coordinate $u$ on $H$).
An important consequence  of the theory identifies the weight filtration for
this MHS with the filtration given by the nilpotent operator
$\log ( T)$
where $$T:H^n(X_t) \ra H^n(X_t)$$ is the positive generator for
the monodromy action of $\pi_1(\D^*)$
induced by the transformation
$u \mapsto u+1$ of the upper half plane.
Furthermore,   $\log (T)$ shifts weights by $-2$ and Hodge types by
$(-1,-1)$. More precisely, we have a map of (mixed) Hodge structures
$$N:=-(1/2\pi i)\log (T): H^n(X_t)(1) \ra H^n (X_t)$$
Using duality, we can therefore view $N$ as an element
in $$(H^n)^*\otimes H^n(X_t) (-1)\simeq H^{2d-n}(X_t)\otimes H^n(X_t)(d-1)$$
generating a trivial substructure. 
Summing over all $n$ and using the K\"unneth formula, 
we get an element
$$[N] \in H^{2d}(X_t \times X_t) (d-1)$$
generating a trivial substructure of the limit mixed Hodge structure on the
cohomology of the product space
$X_t\times X_t$. Now, weight considerations show that
$$[N] \in [W_{2d-2}H^{2d}(X_t \times X_t)] (d-1)$$
thus giving rise to a trivial substructure of
$$[Gr^W_{2d-2}H^{2d}(X_t \times X_t)] (d-1)$$
The situation  here is obviously different from the usual one
for the Hodge conjecture in that we are not viewing the cohomology groups
with the usual pure Hodge structure. Nonetheless, the construction of
Steenbrink's spectral sequence provides a natural context for
interpreting this situation. For this, let $$Z\ra \D$$
be a regular normal-crossing model for $X\times_{\D} X$,
and let $S$ be its special fiber. It will be made up
of smooth components  that cross normally.
Denote by $$S^{[i]}$$ the disjoint union of all $i$-fold
intersections of these components.

Then Steenbrink's spectral sequence expresses
the pure graded pieces of $H^{2d}(X_t \times X_t)$ as made up
of the $E_2$ term of a spectral sequence, whose $E_1$ term consists
of sums of  cohomology groups of the $S^{[i]}$. In particular,
the piece of interest
$$[Gr^W_{2d-2}H^{2d}(X_t \times X_t)] $$
can be expressed as a subquotient of the cohomology of
smooth projective varieties.

We can be more precise. For this, we note  that 
$N$ viewed as a map $$H(X) (1) \ra H(X)$$ commutes with the monodromy action
on the source and target, and hence, the
class $[N]$ is invariant for the monodromy operator on
the (twisted) cohomology of the product. Now, an important consequence of
the theory (the invariant cycle theorem) expresses the monodromy
invariant part of a limit MHS as the image of the MHS of the special
fiber. That is, the natural map
$$H^{2d}(S) \ra H^{2d} (X_t\times X_t)$$
surjects onto the invariant part. On the other hand, since
the morphism of Hodge structures is strict, we get thereby a surjection
$$Gr^W_{2d-2}H^{2d}(S) \ra Gr^W_{2d-2}H^{2d} (X_t\times X_t)^I$$
where we have denoted with the superscript $I$ the monodromy invariant part.
We have
$$
Gr^W_{2d-2}H^{2d}(S) \simeq 
 Ker [H^{2d-2}(S^{[3]})\ra H^{2d-2}(S^{[4]})]/ 
Im[H^{2d-2}(S^{[2]})\ra H^{2d-2}(S^{[3]})]$$
where the maps are alternating sums of restriction maps (with signs determined
by fixing an ordering of the components), weighted suitably
by the multiplicities of the components. Therefore,
there exists a trivial substructure of $H^{2d-2}(S^{[3]})(d-1)$
which maps to the subspace generated by $[N]$, and the Hodge conjecture
would imply that this space can be represented by an algebraic cycle.
The problem of identifying this algebraic cycle
 was posed by Kazuya Kato. Stated more precisely, he pointed out
the following consequence of the Hodge conjecture:

\begin{conj}
There exists an algebraic cycle on $S^{[3]}$ whose cycle class lies in the
kernel of the restriction map to $S^{[4]}$ and which maps
to the class $[N]$ of the monodromy operator discussed above.
\end{conj}

This conjecture 
was proved for the case of curves with double
point degenerations and surfaces with triple point degenerations
in \cite{con} for a specific model $Z$.
In fact, one of the results of the present paper
is that
\begin{thm} The conjecture is independent of the model
$Z$ of $X\times_{\D} X$. That is, if the
conjecture is verified for one normal-crossing model $Z$ then the
conjecture holds for
any other normal-crossing model $Z'$.
\end{thm}

Furthermore, we will prove the conjecture in the very simple case of
semi-stable fibrations arising from isolated quadratic singularities.
That is, assume that $V'\ra \D$ is a projective flat map 
where the singularities of the map are all in the special fiber
$W'$ and are isolated
non-degenerate quadratic. Let $V$ be obtained by blowing up
the singular points of $W'$. Finally, let $X$ be obtained by normalizing
the base change of $V'$ via the squaring map $\D \ra \D$.
Thus, $X$ is semi-stable, and we will have a monodromy class
$[N] \in H^{2d}(X_t\times X_t)(d-1)$.

\begin{thm} Kato's conjecture is true in this setting.
That is, there is an algebraic cycle class in
$H^{2d-2}(S^{[3]})(d-1)$ mapping via the Steenbrink
spectral sequence to the class $[N] $.
\end{thm}

The proof is a simple application of the Picard-Lefschetz formula
and can easily be guessed at by experts. On the other hand,
since the method suggests an approach to the general problem with
possibly considerable
technical advantages over the approach of \cite{con},
we considered it worth writing down explicitly.

This problem has an analogue in the case of mixed-characteristic
degenerations that we do not treat here. But we note that
the monodromy operator is the subject of one of the most
important open problems of arithmetic geometry, namely, the
weight-monodromy conjecture. 
It is our hope (perhaps futile) that the
algebraicity of Kato's conjecture may eventually contribute
to its resolution.

\section{Proofs}
{\em Proof of theorem 1.}

 By the weak factorization theorem \cite{Wl},
any two normal crossing models (that is,
regular models having normal crossings in the special fiber)
can be connected by a sequence of diagrams
$$\begin{array}{ccccc}
 & & Z' & & \\
 & \swarrow & & \searrow & \\
Z& & & & Z''
\end{array}$$
where each arrow is a blow-up along a non-singular
center in the special fiber. Also, we can assume
all the models occuring in the
diagram are normal crossing and that the
center of each blow-up  meets the components of the
special fiber normally.
Thus, by iterating along the diagram,
we need only show  that the algebraicity of
the monodromy class $[N]$ for the two models
$Z$ and $Z'$ are equivalent.
Denote the map by $f:Z'\ra Z$ and let $E$ be the
exceptional divisor. 

Let us choose some complex computing the cohomology
of the Steenbrink complexes for the two models,
say that constructed using $C^{\infty}$ differential forms.
That is, we let $A$ and $A'$ be the global $C^{\infty}$
versions of the Steenbrink complexes for $Z$ and $Z'$,
where the sheaf of log differential forms is replaced
by global $C^{\infty}$ differential forms with log poles.
We clearly have a map $f^*:A \ra A'$ which preserves the
weight and Hodge filtrations.
Thus, we have a map of the Steenbrink spectral sequences
associated to the two models.
\[
f^*: E_1^{-r,q+r}(Z) \to E_1^{-r,q+r}(Z'),\quad f_!: E_1^{-r,q+r}(Z')
 \to E_1^{-r,q+r}(Z).
\]
By the description of the associated graded pieces via
residues it follows that the pull-back map at the
$E_1$ level is given by usual pullback in cohomology
for the components corresponding to various components
of strata.
That is, at the rational level, these morphisms are described by a direct sum 
of maps as 
\begin{align*}
f^*&: H^{q-r-2k}(S^{[2k+r+1]},\Q(-r-k)) \to 
H^{q-r-2k}({S'}^{[2k+r+1]},\Q(-r-k)),\\
\end{align*}
where $k$ is an integer satisfying $k \ge\text{max}(0,-r)$. 

Since both spectral sequences are associated to a mixed
Hodge complex converging to the cohomology of
the generic fiber, we see that $f^*$ is a quasi-isomorphism
when we regard $E_1(Z)$ and $E_1(Z')$ as complexes.
We can in fact construct a splitting $f_{!}:E_1(Z') \ra E_1(Z)$
as follows.
The components of $(S')^{[i]}$ are of the form
$$S_{a_{1}}'\cap \cdots \cap S_{a_i}'=(S_{a_{1}}\cap \cdots \cap S_{a_i})'$$ 
and
$$ S_{a_{1}}'\cap \cdots \cap S_{a_{i-1}}'\cap E=
(S_{a_{1}}\cap \cdots \cap S_{a_{i-1}})'\cap E$$
where $S_k'$ is the strict transform of the
component $S_k$ of $S$. Here, we have chosen the ordering of the components
of $S'$ in such a way that it is compatible with the ordering of
the $S_k$'s and so that $E$ is the last element.
Then we define
$f_{!}$ on 
$$H^n(S_{a_{1}}'\cap \cdots \cap S_{a_i}')$$
as the usual push-forward to
$$H^n(S_{a_{1}}\cap \cdots \cap S_{a_i})$$
and on
$$H^n(S_{a_{1}}'\cap \cdots \cap S_{a_{i-1}}'\cap E)$$
as zero.

Let us check that $f_{!}$ commutes with the differentials on
$E_1(Z')$ and $E_1(Z)$. On both complexes the differential
$d$ has two
components $d_1$ and $d_2$, the first consisting of restriction maps
and the other of Gysin maps. (There is a slight complication
coming from the weighting for the multiplicities. But these
will be the compatible on the two complexes because the
center of the blow-up meets the divisor $S$ normally, so
we will ignore them.)

Assume we start from a class 
$$c\in H^n(S_{a_{1}}'\cap \cdots \cap S_{a_i}')$$
The various restriction maps will have components in groups of the form
$$H^n(S_{a_{1}}'\cap \cdots \cap S_{a_i}'\cap S_{a_{i+1}})'$$
for which we have a commutative diagram
$$\begin{array}{ccc}
H^n(S_{a_{1}}'\cap \cdots \cap S_{a_i}') &\ra &
H^n(S_{a_{1}}'\cap \cdots \cap S_{a_i}'\cap S_{a_{i+1}}')\\
\downarrow & & \downarrow \\
H^n(S_{a_{1}}\cap \cdots \cap S_{a_i}) &\ra &
H^n(S_{a_{1}}\cap \cdots \cap S_{a_i}'\cap S_{a_{i+1}})
\end{array}$$
while the restriction to components of the form
$$H^n(S_{a_{1}}'\cap \cdots \cap S_{a_i}'\cap S_{a_i}\cap E)$$
composes to zero under $f_{!}$. This shows that
$f_{!}(d_1(c))=d_1f_{!}(c)$.
The other compatibility,
$f_{!}(d_2(c))=d_2f_{!}(c)$
just follows from the functoriality of push-forwards.

Now start from
$$c\in H^n(S_{a_{1}}'\cap \cdots \cap S_{a_{i-1}}'\cap E).$$
 Then both $f_{!}d_1$ and $d_1f_{!}$ are zero. 
Also, $d_2f_{!}(c)=0$ by definition. So to conclude, we need
only check that $f_{!}d_2(c)$ is zero.
This is clear for the components of $d_2$ corresponding to
Gysin maps of the form
$$H^n(S_{a_{1}}'\cap \cdots \cap S_{a_{i}}'\cap E)
\ra H^{n+2}(S_{b_{1}}'\cap \cdots \cap S_{b_{i-1}}'\cap E)$$
For the map
$$H^n(S_{a_{1}}'\cap \cdots \cap S_{a_{i}}' \cap E)
\ra H^n(S_{a_{1}}'\cap \cdots \cap S_{a_i}')$$
note that 
\[
H^{n+2}(S_{a_{1}}'\cap \cdots \cap S_{a_i}')\simeq H^{n+2}(S_{a_{1}}'\cap \cdots \cap S_{a_i}')\oplus
H^n(S_{a_{1}}'\cap \cdots \cap S_{a_i}'\cap E)
\]
by the formula for the cohomology of a blow-up
and that $f_!$ is just projection onto the first component
while the second component is the image of the Gysin map.
Thus, $f_!d_2(c)=0$ and we are done.

We can write now
$$E_1(Z')\simeq E_1(Z)\oplus K$$
as complexes where $K$ is the kernel of $f_!$ and is acyclic.
In particular, we see that both $f^*$ and $f_!$
induce maps between
$$Ker (d:E_1(Z) \ra E_1(Z))$$
and $$Ker (d:E_1(Z') \ra E_1(Z'))$$ which are compatible with
with maps from either to $$Gr^W(H^*(X_t\times X_t)).$$
Finally, noting that both $f^*$ and $f_!$
take algebraic classes to
algebraic classes allows us to conclude the proof.
\medskip

{\em Proof of theorem 2.}

We start by reviewing the relevant portions of the Picard-Lefschetz
formalism (\cite{De3} XV).

To restate our assumptions,
let $V'\ra \D$ be a proper flat map, smooth over $\D^*$. We now assume
that $V'$ has odd fiber dimension $d$ and that the singularities
of the special fiber are all quadratic and isolated,
coming from a local expression $\underline{x}\ra Q(\underline{x})$ for our
map $V' \ra \D$, where $Q$ is a non-degenerate
quadratic form. We note
at this
point that the case of even fiber dimensions is trivial
with respect to our problem
since the monodromy operator is then zero.
Let $V\ra V'$
be the blow-up of $V'$ along the singularities of the special
fiber.Denote by $W'$ and $W$ the special fibers of
$V'$ and $V$ respectively.
 We see that $W$ is
normal crossing with
$s+m$ components, where $m$ is the number of components of
$W'$ and $s$ is the number of singularities in
the special fiber $W'$. Note that $m=1$ unless $d=1$ with our assumptions.
For the rest of this paper, we will therefore assume that
$d>1$ and hence $m=1$, since the arguments for
 $d=1$ simply involve omitting some steps
in an obvious manner, and another proof is given in \cite{con}
in any case.

Let 
$X\ra \D$ be obtained from $V$  
by changing 
base with respect to the map
$\D\ra \D, \ \ t \mapsto t^2$ and normalizing. Thus, 
$X$ is a smooth manifold and the special fiber $Y$ is
reduced normal crossing. The set of its components
are in bijection with the components of
$W$, and we label them $\{Y_i\}$, $\{W_i\}$,
for some label set to be determined later,
so that the map $X\ra V$ maps $Y_i$ to $W_i$.
Let $T'$ be the monodromy operator
on the cohomology of $V'_t$, any smooth fiber of $V'$. 
Thus, $T'$ acts non-trivially
only on $H^d(V'_t)$. Then the Picard-Lefschetz formula
says
$$T'(x)=x-\Sigma_p(\delta_p,x)\d_p$$
where $p$ runs over the singularities of $V'_0$ and
$\d_p \in H^d(V'_t)$ is the vanishing cycle associated to $p$.
$\d_p$ is determined (up to sign) as follows:
we have the composed map
$$H^d_p(W') \ra H^d(W') \ra H^d(V'_t)$$
and $H^d_p(W')$ is free of rank one. Take $\d_p$
to be the image of any generator of $H^d_p(W')$.
The formula is clearly independent of the choice of
generator. Now, if $T$ denotes the monodromy operator for
$X$ acting on $H^d(X_t)=H^d(V'_{t^2})$, 
then we have $T=(T')^2$ and $T'$ acts trivially
on the vanishing cycles, so we get
$$T(x)=x-\Sigma_p 2(\delta_p,x)\d_p$$
Therefore, the log of the monodromy is given by
$$\log T (x)=T(x)-x=-\Sigma_p 2(\delta_p,x)\d_p$$
since $(T-1)^2=0$, and the cohomology class $[N]$ in
$$H^{2d}(X_t\times X_t)(d-1)$$ defined
by $N=-(1/2\pi i)\log T$ is the class
$$(2\pi i)^{d-1}
2\Sigma_p \d_p \otimes \d_p \in H^d(X_t) \otimes H^d(X_t) (d-1) \subset
H^{2d}(X_t\times X_t)(d-1)$$

We will show, in fact, that every class $(2\pi i)^{d-1} \d_p\otimes \d_p$
is `algebraic' in the sense described in the
introduction.
For this, we take  $Z$ to be the blow-up of
$X\times_{\D} X$ along its singular locus,
which we will describe more precisely below.
We have a commutative diagram
$$\begin{array}{ccc}
X_t\times X_t &\hra &Z\\
|| & &\da \\
X_t \times X_t& \hra &X\times_{\D} X
\end{array} $$
inducing a commutative diagram of
maps in cohomology
$$\begin{array}{ccccc}
H^i(X_t \times X_t) &\leftarrow & H^i(Z) &\simeq &H^i(S)\\
|| & &\uparrow & & \uparrow\\
H^i(X_t \times X_t )& \leftarrow & H^i(X\times_{\D} X) &\simeq & H^i(Y\times Y)
\end{array} $$
and we need to show that the class $[N]$ in
$H^{2d} (X_t\times X_t)(d-1)$ comes from an algebraic
class in $H^{2d-2}(S^{[3]})(d-1)$.
To get this, we need to analyze a bit
the original Lefschetz fibration
$V'$. The class $\d_p \in H^d(V'_t)$ comes from
a class  $\e_p \in H^d(W')$
which in turn comes from the (rank 1) local cohomology
$H^d_p(W')$ via the map
$H^d_p(W')\ra H^d(W')$. On the blowup $V$,
the special fiber $W$ contains a component $W_p$
mapping to the point $p$ for each $p$ and
there is a component $W_0$ given by
the strict transform of $W'$. Also,
$W_{0p}:=W_0\cap W_p$ is the exceptional divisor of the (blow-up) map
$W_0 \ra W'$.
There is a surjection
$H^{d-1}(W_{0p}) \ra  H^d_p(W')$
defined as follows:
Let $U_p$ be a contractible neighborhood of $p$ in $W'$
and $N_p$ its inverse image in $W_0$. Then $N_p$ is
homotopy equivalent to $W_{0p}$ while
$N_p-W_{0p} \simeq U_p-p$. So we have maps
$$H^{d-1}(W_{0p}) \simeq H^{d-1}(N_p) \ra H^{d-1}(N_p-W_{0p}) \simeq
H^{d-1}(U_p-p) \ra H^d_p(U_p)=H^d_p(W')$$
But $H^{d-1}(W_{0p})$ is the middle degree cohomology of the
smooth quadric $W_{0p}$ and therefore (\cite{De3} XII.3.3)
is generated
by a power of the hyperplane class and one other class $\a'_p$
(corresponding to the maximal isotropic
subspaces for the quadratic form defining $W_{0p}$)
representing the primitive quotient. We can take
$e'_p\in H^d(W')$ to be the image of $\a'_p$.
Now, going back up to
$X$, we therefore get classes
$\a_p \in H^{d-1}(Y_{0p})$
where $Y_{0p}$ is the intersection
$Y_0\cap Y_p$ mapping to
the pullback $e_p$ of $e'_p$ to $H^d(Y)$.
By construction,
$2\Sigma_p e_p\otimes e_p \in H^{2d}(Y\times Y)$
maps to $[N]$ via the map to the generic
fiber of $X\times_{\D } X$. Denote by $c_p$ the class
in $H^{2d}(S)$ obtained by pulling back
$e_p\otimes e_p$ which, therfore, is the same as
$$p_1^*(e_p)\cup p_2^*(e_p) \in H^{2d} (S)$$
where $p_i$ are the composition of the blow-up
map $S\ra Y\times Y$ with the projection maps to the two factors.
Thus, $(2\pi i)^{d-1}2\Sigma_p c_p$ is now the class in the cohomology of
$S$ which maps to the class $[N]$.
We wish to show that $c_p$ comes from an algebraic class
in $H^{2d-2}(S^{[3]})(d-1)$.
For this, we need to analyze a bit the geometry of the blow-up
$Z\ra X\times_{\D} X$. For each index $i$, either $0$ or $p$,
we have a component $Y_i$ of $Y$, giving us components
$Y_i\times Y_j$ of $Y\times Y$. We also have the singular
locus
$$Y_{ij} \times Y_{lk}$$
(using the obvious notation of double indices for
the double intersections) of $X\times_{\D} X$ which we need to blow up to
get $Z$. We will find the appropriate classes in the
inverse image of
$Y_{0p} \times Y_{0p}$.

We need only deal with each $p$ separately, so
we will now fix a $p$ of interest and
index the components of $S$ as follows:
$S_0$ refers to the exceptional divisor mapping
to $Y_{0p} \times Y_{0p}$. $S_1, S_2, S_3, S_4$
are then the strict transforms of $Y_0\times Y_0,
Y_0 \times Y_p, Y_p \times Y_0, $ and $Y_p \times Y_p$
respectively. If we note that the center of
the blow-up is a Cartier divisor on 
$(Y_0\times Y_0 )\cap (Y_0 \times Y_p)$, we get that
$S_{012}=S_0 \cap S_1 \cap S_2$ maps isomorphically
to the center $Y_{0p} \times Y_{0p}$. Similarly,
$S_{013}$,$S_{024}$, and $S_{034}$ map isomorphically
to  $Y_{0p} \times Y_{0p}$. Furthermore, these are the only
non-empty triple intersections.
Now we claim that there is an algebraic class in
$$  H^{2d-2}(S_{013})(d-1)$$
which maps to $(2\pi i)^{d-1}c_p$. To see this, we need some
explicit computation with cocycles which describe the
various maps involved.

First, let us describe the map 
$$f_1: H^{d-1}(W_{0p}) \ra H^d (W)$$
obtained as the composite
$$H^{d-1}(W_{0p}) \ra H^d(W')\ra H^d(W)$$
We claim it is the same as the map $$f_2:H^{d-1}(W_{0p}) \ra H^d (W)$$
obtained from the spectral sequence for the
hypercover
$$\coprod W_i \ra W$$
To see this, we recall how to compute  $f_1$
at the level of cocyles. Denote by $f$ the  map
$W\ra W'$ as well as its restriction to various subsets.
Given a class $[x]\in H^{d-1}(W_{0p})$, let $x$ be a cocycle representing it.
We can lift it to a cocycle $y$ on $N_p$ and then restrict it
to $y_0$ on $N_p-W_{0p}$ which can then be expressed as
$f^*(y_0')$ for some cocycle $y_0'$ on $U_p-p$. To get the image in
$H^d_p(U_p)$, we just consider the relative cocycle
$(0, y_0')$ in the cone complex computing local cohomology.
Now we need to find a cocycle representing the
class in $H^d_p(W')\simeq H^d_p(U_p)$. To do this, we solve
$[(0,y_0')]=[(z,w)]|(U_p,U_p-p)$ where $z$ is a $d$-cocycle on $W'$
and $w$ is a $d-1$ cochain on $W'-p$ such that 
$z|(W'-p)=dw$ and the equality of classes
means that there is a pair $(a,b)$ consisting of a $d-1$ cochain
$a$ on $U_p$ and a $d-2$ cochain $b$ on $U_p-p$
solving
$$z|U_p=da, \ \ w|(U_p-p)=y_0'+a|(U_p-p)-db|(U_p-p)$$
Then $f_1([x])$ is the class of $f^*(z)$  on $W$.

In the above, notice that we can take
$b=0$. This is done as follows: Consider the covering
$W'=(W'-p)\cup (U_p)$. Then the complex underlying the Mayer-Vietoris
sequence for this covering allows us to write
$y_0'=w-a$ for a cochain $w$ that comes from $W'-p$ and
$a$ that comes from $U_p$. Since, $y_0'$ is a cocycle,
we get $dw=da$ on $U_p-p$, giving us a class $z$ on $W'$
satisfying $z=dw$ on $W'-p$ and $z=da $ on $U_p$.
Then $(z,w)$ is a cocycle for relative cohomology
satisfying the constraints above for $b=0$.

Pulling back to $W_0$ using $f$,
we get an equality 
$$f^*(w)|(N_p-W_{0p})=(y+f^*(a))|(N_p-W_{0p})$$
the point being that on $W_0$, $f^*(y_0')=y_0$
extends to the class $y$ on $N_p$.
So we get a class $v$ on $W_0$, defined to be
$f^*(w)$ on $W_0-W_{0p}$ and $y+f^*(a)$ on $N_p$, with the property that
$v|W_{0p}=x$. Then $f_1([x])$ is represented by the cocycle
$f^*(z)$ on $W$, which is characterized by the
property of restricting to $dv$ on $W_0$
and 0 on $W_p$.

On the other hand, we compute $f_2$ as follows.
Let $w_{0p}$ be a $d-1$ cocycle on $W_{0p}$. Find cochains
$w_0$ and $w_p$ on $W_0$ and $W_p$ such that
$w_{0p}=w_0-w_p$ on $W_{0p}$. So we have a cochain
$(w_0,w_p)$ in the direct sum of the complexes computing the
cohomology of $W_0$ and $W_1$.
Then compute its
boundary to get $(dw_0,dw_1)$ which comes from a $d$ cocycle
$w $ on $W$, that is, $w|W_0=dw_0$ and $w|W_p=dw_p$.
Then $f_2([w_{0p}])=[w]$.
To compare the two maps, start with the cocycle
$x$ on $W_{0p}$. Above, we constructed a cochain
$v$ on $W_0$ such that $v|W_{0p}=x$.
So we get the element $(v,0)\in C^{d-1}(W_0)\oplus C^{d-1}(W_p)$
such that $v-0=x$ on $W_{0p}$. Take the differential
to get $(dv,0)$.  As noted above, $f_1([x])$ is represented by the
cocycle
$f^*(z)$ which is a class on $W$ that restricts to zero on
$W_p$ and $dv$ on $W_0$. This shows that $f_1([x])=f_2([x])$.

Next, consider the map
$$H^{d-1}(Y_{0p})\otimes H^{d-1}(Y_{0p}) \ra H^{2d} (Y\times Y) \ra H^{2d}(S)$$
In order to analyze the class $c_p$,
we need to compare this with the map
$$ H^{2d-2}(S_{012})\oplus H^{2d-2}(S_{013})\oplus
H^{2d-2}(S_{024})\oplus H^{2d-2}(S_{034})\ra H^{2d}(S)$$
coming from the spectral sequence of the hypercover
$$\coprod S_i \ra S$$
(Here, since we're speaking of a full hypercover of $S$
although we will be interested only in the part coming from
$S_0,S_1,S_2,S_3,S_4$, we assume that we have extended the labelling
to all of the components of $S$ in some manner.)
Take the class $c_p=p_1^*(e_p)\otimes p_2^*(e_p) \in H^{2d}(S)$
and restrict to the components $S_0, S_1,S_2,S_3,S_4$.
For convenience of notation, we will suppress  the pull-back
maps and write, for example, $c_p=e_p e_p$.
Recall that $S_0$ is the exceptional divisor
of the blow-up $Z\ra X\otimes_{\D}X$ while
$S_1,S_2,S_3,S_4$ are the strict transforms of
$Y_0\times Y_0,
Y_0 \times Y_p, Y_p \times Y_0, $ and $Y_p \times Y_p$
respectively. 
The center of the blow-up is
$Y_{0p}\times Y_{0p}$. Furthermore, we have represented
$e_p$ as a cocycle characterized by
  $e_p|Y_0=dv$, $e_p|Y_p=0$, where
$v$ is a cochain on $Y_0$ satisfying $v|Y_{0p}=x$.
Now, the restriction of $c_p=e_pe_p$ to $S_0$ is the same
as taking $p_1^*(e_p)p_2^*(e_p)$ on $Y\times Y$, restricting to
$Y_{0p}\times Y_{0p}$ and pulling back to the exceptional divisor. So
$c_p|S_0=0$.
Arguing this way with the various restrictions and using the fact that
$dv|Y_p=0$, 
we see that the only component that
survives is $c_p|S_1=dvdv=d(vdv).$ That is,
we have the class
$$(0,vdv,0,0,0)\in C^{2d-1}(S_0)\oplus C^{2d-1}(S_1)\oplus 
C^{2d-1}(S_2)\oplus C^{2d-1}(S_3)\oplus C^{2d-1}(S_4)$$ 
whose differential is
the restriction of $c_p$ to the components.
Now,taking the Cech differential, by the same argument as above
using  $dv|Y_p=0$, we get that the only component
among the double intersections that survives is
$-vdv$ on $S_{13}$. This is the differential of
the $2d-2$ cochain $-vv$ on $S_{13}$.
The Cech differential of this element has
non-zero component only on $S_{013}$ which is equal to
$-\pi^*(p_1^*(\a_p)p_2^*(\a_p))$
for the isomorphism $\pi:S_{013} \ra Y_{0p}\times Y_{0p}$.
We get that $(2\pi i)^{d-1}$ times this
class is clearly algebraic, since $(2\pi i)^{(d-1)/2}a_p$ is algebraic.
Then  $(2\pi i)^{d-1}c_p$, which is the image of
$$(2\pi i)^{d-1}\pi^*(p_1^*(\a_p)p_2^*(\a_p))$$
 under the map $$H^{2d-2}(S_{013})(d-1) \ra H^{2d}(S)(d-1)$$
is algebraic, and hence, so is $[N]$.
We remark that the sign may differ from that given above
depending on a convention for ordering
various components.
\medskip

{\em Acknowledgements:} M.K. was supported in part by the
NSF and  C.C. was supported in part by
 NSERC grant 72016789.
M.K. is grateful to the Korea Institute for Advanced Study
and the University of Toronto for their hospitality.
Both authors wish to thank the MPIM where part of this
work was completed.

{\footnotesize C.C.: DEPARTMENT OF MATHEMATICS, UNIVERSITY
OF TORONTO,
 TORONTO, ONTARIO, CANADA M5S 3G3, EMAIL: kc@math.toronto.edu}

{\footnotesize M.K.: DEPARTMENT OF MATHEMATICS, UNIVERSITY
OF ARIZONA, TUCSON, ARIZONA 85721, U.S.A, EMAIL: kim@math.arizona.edu}
\end{document}